\documentclass[letterpaper, 10 pt, conference]{ieeeconf}  

\IEEEoverridecommandlockouts                              

\usepackage{amsmath,amsfonts}
\usepackage{hyperref} 
\usepackage{algorithm}
\usepackage{algpseudocode}
\usepackage{graphicx}
\newcommand{\cmmnt}[1]{}

\algnewcommand\algorithmicforeach{\textbf{for each}}
\algdef{S}[FOR]{ForEach}[1]{\algorithmicforeach\ #1\ \algorithmicdo}
\renewcommand{\P}{\mathcal{P}}
\newcommand{\T}{\mathcal{T}}
\DeclareMathOperator*{\argmin}{arg\,min}

\newcommand{\bpi}{\boldsymbol{\pi}}
\newcommand{\bp}{\boldsymbol{p}}
\newcommand{\bG}{\boldsymbol{G}}
\newcommand{\bh}{\boldsymbol{h}}

\title{\LARGE \bf
A Stackelberg Game Approach to Control\\ the Overall Load Consumption of a Residential Neighborhood
}

\author{Erhan Can Ozcan$^{1}$ and Ioannis Ch. Paschalidis$^{2}$
\thanks{The research was partially supported by the DOE under grants DE\-EE0009696 and DE\-AR\-0001282, the NSF under grants CCF\-2200052, DMS\-1664644, and IIS\-1914792, and by the ONR under grants N00014\-19\-1\-2571 and N00014\-21\-1\-2844. }
\thanks{This work has been submitted to the IEEE for possible publication. Copyright may be transferred without notice, after which this version may no longer be accessible.}
\thanks{$^{1}$Erhan Can Ozcan is with Division of Systems Engineering, Boston University, Boston, MA 02215, USA
        {\tt\small cozcan@bu.edu.}}%
\thanks{$^{2}$Ioannis Ch. Paschalidis is with the Department of Electrical and Computer Eng., Division of Systems Eng., and Dept. of Biomedical Eng., Boston University, 8 St. Mary’s St., Boston, MA 02215, USA
        {\tt\small yannisp@bu.edu.}}%
}

\begin{document}

\maketitle
\thispagestyle{empty}
\pagestyle{empty}

\begin{abstract}

This paper formulates a Stackelberg game between a coordination agent and participating homes to control the overall load consumption of a residential neighborhood. Each home optimizes a comfort-cost trade off to determine a load schedule of its available appliances in response to a price vector set by the coordination agent.\cmmnt{This paper formulates a Stackelberg game between a coordination agent and participating homes to control the overall load consumption of a residential neighborhood, where each home determines the optimal load schedule of its available appliances in response to the price vector announced by the coordination agent according to its own comfort related preferences.} The goal of the coordination agent is to find a price vector that will keep the overall load consumption of the neighborhood around some target value. After transforming the bilevel optimization problem into a single level optimization problem by using Karush-Kuhn-Tucker (KKT) conditions, we model how each home reacts to any change in the price vector by using the implicit function theorem. By using this information, we develop a distributed optimization framework based on gradient descent to attain a better price vector. We verify the load shaping capacity and the computational performance of the proposed optimization framework in a simulated environment establishing significant benefits over solving the centralized problem using commercial solvers.
\end{abstract}

\section{INTRODUCTION}

Maintaining a balance between electricity demand and supply is extremely important to ensure the reliability of power grids, and the main challenge utility companies are facing is to meet the electricity demand during peak hours \cite{my_ref1}. In order to overcome this issue, utility companies dispatch additional power plants with high operating costs during peak hours\cmmnt{, which causes a huge financial loss}. The increasing production from intermittent renewable sources augments the need for such reserve capacity. Moreover, as these power plants mostly rely on fossil fuels to generate electricity, they increase carbon dioxide emissions and accelerate climate change\cmmnt{significantly, thus damaging the environment irreversibly } \cite{my_ref2}. According to the report published by the U.S. Energy Information Administration, the energy consumption in residential dwellings will increase by $22\%$ until 2050, and the main underlying factor behind this increase is the population growth \cite{my_ref3}.\cmmnt{According to the report published by the U.S. Energy Information Administration, the energy consumption in residential dwellings will increase by $22\%$ until 2050, and this increase is mainly driven by the population growth and greater reliance on electronic appliances \cite{my_ref3}.} As the increase in the number of end users makes the electricity peak load problem more severe \cite{my_ref1}, the studies focusing on controlling the load consumption in residential homes and buildings is of utmost importance for the society.

Deployment of smart meters and {\em{Home Energy Management Systems}} (HEMS) to manage home energy consumption and energy costs has become popular after the developments in information and communication technologies \cite{my_ref4}. The change in the energy usage of a consumer in response to varying electricity prices or to incentive payments is referred as {\em{Demand Response}} (DR) \cite{my_ref5}. \cmmnt{and various DR programs have been introduced in \cite{my_ref5}.}

This paper proposes a demand response program with dynamic pricing, which has shown promising results to decrease the peak load in real-world applications\cmmnt{ by shaping the load profiles of individual homes} \cite{my_ref6}, \cite{my_ref7}.
\cmmnt{The demand response programs with dynamic pricing have shown promising results to decrease the peak load by shaping the load profiles of individual homes \cite{my_ref6}, \cite{my_ref7}.}
A Stackelberg game formulation between the energy service provider and homes is introduced in \cite{my_ref8}, where the user comfort is represented by a utility function. However, it is not possible to develop a single utility function to guarantee user preferences are always met by the resulting load plan. Since the participation of users is a key factor in the success of demand response programs, methods failing to maintain the user comfort may not realize their full potential \cite{my_ref9}. Hence, each appliance in HEMS should be modeled with enough details so that the user specific preferences can be reflected as closely as possible. In \cite{my_ref10}, the authors propose a Stackelberg game maximizing the profit of the utility company and minimizing the electricity bill of the homes. Although different appliances such as shiftable, non-shiftable and curtailable are considered, these appliances do not have enough details to reflect the real dynamics. In \cite{ref_ACEEE}, the authors propose a Stackelberg game where the electricity retailer wants to control the overall load consumption of the neighborhood with sparse deviations from a target value. However, the resulting optimization problem is a non-convex {\em{Quadratically Constrained Quadratic Programming}} (QCQP) problem, which can be computationally challenging to solve, especially in large neighborhoods. Therefore, employing this model in real-world applications may not be practical unless efficient computational methods are developed to solve it. In \cite{prev_paper}, a {\em{Mixed Integer Linear Programming}} (MILP) is proposed to control the load consumption of the neighborhood without considering the electricity price. In this scenario, the proposed formulation is not a Stackelberg game but a single level optimization problem, and the authors solve it by Dantzig-Wolfe decomposition. However, the application of this strategy in real life may require participating homes to sign long term contracts with fixed rates, which makes it impossible to take advantage of the low electricity price when there is ample supply in the electricity market.

In this study, we formulate a bilevel optimization problem where the coordination agent announces a price vector, and then each home determines the optimal load schedule of the appliances in response to that price vector according to its own comfort-related preferences. Since each home solves a quadratic programming problem to find its optimal load schedule, we transform the bilevel optimization problem into a single level optimization problem by using Karush-Kuhn-Tucker (KKT) optimality conditions. Solving the reformulated single level problem by using commercial solvers is the state-of-the-art strategy to solve bilevel optimization problems \cite{my_ref12}, \cite{my_ref13}. However, we stress that this solution strategy has two major drawbacks preventing it from being used in practice.\cmmnt{Although commercial solvers can deal with the reformulated problem by using {\em{Mixed Integer Quadratic Programming}} (MIQP) techniques, we stress that this solution strategy has two major drawbacks preventing it from being used in practice.} First, the problem complexity grows exponentially as the target community becomes larger, which makes attaining an effective solution in a reasonable amount of time difficult. Second, the coordination agent needs to know both the preferences and the sensitive personal data of homes to be able to formulate the aforementioned optimization problem, which can raise some data privacy concerns among users. \cmmnt{Second, the proposed solution technique requires participating homes to share their preferences and sensitive personal data with the coordination agent at each time interval, which can raise some data privacy concerns among users.} In order to address the problems related to efficiency and data privacy, we develop a gradient based distributed optimization framework similar to the approach introduced in \cite{my_ref12}.

We compare the proposed optimization framework with a commercial solver in terms of the optimization time and the solution quality. Our experiments demonstrate that the commercial solver fails to solve the centralized problem within the desired time interval when the optimization problem includes more than a few homes. In most of the cases, it cannot even find an initial feasible solution. On the other hand, the distributed optimization framework is able to provide an effective solution in a short period of time. Furthermore, we show that better solutions can be obtained by utilizing mini-batch gradient descent optimization techniques.

Our contributions can be summarized as follows: $1)$ We formulate a Stackelberg game between a coordination agent and participating homes to minimize the peaks observed in the power profile of the community without severely affecting consumer comfort. $2)$ We develop a gradient-based distributed optimization approach to solve the proposed problem efficiently, and compare it with a commercial solver. $3)$ We show that stochastic unbiased estimators for the gradient exist, and using stochastic gradients can bring benefits while optimizing the formulated problem.  

\textbf{Notation}: We denote vectors by bold lowercase letters (e.g., $\bp$ and $\bpi$), and matrices by bold uppercase letters (e.g., $\boldsymbol{G}$). We use prime \cmmnt{symbol }to denote the transpose operation. For space-saving reasons, we write ${\boldsymbol{x}}=\left(x(1)\ ... \ x(n)\right)$ for a column vector $\boldsymbol{x}$ in ${\mathbb R}^n$.

\section{MODELING HOME APPLIANCES}
\label{sec_app}

The electricity consumption in residential dwellings is largely attributed to space heating and cooling, water heating, electric vehicle (EV) charging and some other routine home appliances. Therefore, controlling these loads can be an effective way to deal with the peak load problem. In order to ensure the user comfort at any time, each home is equipped with an embedded {\em{Model Predictive Controller}} (MPC) that controls the activity of available appliances in accordance with the user preferences over a time horizon. The appliances considered in this study are the heating, ventilation, and air conditioning (HVAC) system, electric water heater, electric vehicle, and ON/OFF type basic appliances such as washing machine, dryer, and oven. We consider appliance models with the same properties used in \cite{prev_paper} except minor simplifications. In this study, the comfort related constraints of the users are defined by using only linear inequalities, and the full set of constraints with detailed explanation can be found in \cite{this_paper_longer}. \cmmnt{\textbf{or Appendix!!}.}\cmmnt{The full constraint set of each appliance is available in Appendix}
\cmmnt{
\textbf{please see the policy of IEEE L-CSS on referring external files at:} \href{http://ieee-cssletters.dei.unipd.it/Page\_authors.php?p=1}{click to me and check external material/files section}. Do we have enough detail about appliances in this section so that we meet the criteria?}

\section{PROPOSED OPTIMIZATION PROBLEM}

Consider a power grid managed by a coordination agent serving $N$ participating homes equipped with $M$ appliances introduced previously. The total load consumption of the community can be controlled by formulating a Stackelberg game where the coordination agent is the leader and the homes are followers. In this game, the coordination agent announces a price vector, then homes automatically adjust the power consumption schedule of the available appliances, minimizing a cost function that accounts for their individual preferences.

Suppose that ${p_{ij}}(t)$ is a decision variable denoting the load consumption of home $i$ for appliance $j$ at time point $t$,\cmmnt{consumed by appliance $j$ in home $i$ at time point $t$.} and \[\boldsymbol{p}_i=\left(p_{i1}(0)\ ...\ p_{i1}(K-1)\ ... \ p_{iM}(0)\ ...\ p_{iM}(K-1)\right)\] is a vector of decision variables denoting the load schedule of home $i$ over the next $K$ time intervals. Each home defines a set of comfort related linear constraints according to its preferences, and these constraints form the polyhedron $\P_i = \{\boldsymbol{p}_i\ |\ \bG_i\boldsymbol{p}_i\leq \bh_i\},\ \forall i$, where $\bG_i$ and $\bh_i$ denote the constraint matrix and the right hand side vector for home $i$, respectively. Hence, the comfort of home $i$ can be ensured as long as $\boldsymbol{p}_i \in \P_i$. Although determining $\boldsymbol{p}_i$ by solving an optimization problem, which minimizes the electricity bill according to the announced price vector, $\boldsymbol{\pi}=\left(\pi(0)\ ...\ \pi(K-1)\right)$, over the set $\P_i$, is possible, there is an oversight in this strategy.\cmmnt{this can be insufficient to reflect the user behaviour.} Each home has different consumption habits, thus a desirable load schedule for each appliance. As any deviation from this schedule causes a social discomfort for the users, these deviations must be taken into account as well as the electricity bill while determining the load schedule of homes. Then, for a given price vector $\boldsymbol{\pi}$, the objective function of home $i$ can be expressed as follows:
\begin{flalign}
\notag f_i(\boldsymbol{p}_i, \boldsymbol{\pi})&=\sum_{j=1}^M\sum_{t=0}^{K-1} p_{ij}(t)\pi(t)\\
&+\sum_{j=1}^M\sum_{t=0}^{K-1}c_{ij}\left(p_{ij}(t)-\overline{p_{ij}}(t)\right)^2,\label{opt_eq_1} 
\end{flalign}
where $c_{ij}$ represents the importance of appliance $j$ for home $i$, and $\overline{p_{ij}}(t)$ represents the desirable load of appliance $j$ in home $i$ at time point $t$. While the first term in Equation \eqref{opt_eq_1} calculates the electricity cost, the second term considers the cost incurred due to deviations from the desirable load curve. Accordingly, home $i$ solves:
\begin{equation}
\label{opt_eq_home}
\boldsymbol{p}_i^*:=\argmin_{\boldsymbol{p}_i\in \P_i}\ f_i(\boldsymbol{p}_i, \boldsymbol{\pi}).
\end{equation}
On the other hand, the coordination agent wants to keep the total load consumption of the community under control while minimizing the total social discomfort of the community by choosing an appropriate electricity price vector $\boldsymbol{\pi}$. Therefore, for a given $\boldsymbol{\pi}$, $\boldsymbol{p}^*$ denotes the optimal load consumption vector of the neighborhood, e.g., $\boldsymbol{p}^*=\left(\boldsymbol{p}_1^*\ ...\ \boldsymbol{p}_N^*\right)$, and the objective function of the coordination agent can be expressed as follows:
\begin{flalign}
\notag f(\boldsymbol{p}^*, \boldsymbol{\pi}) :=&\sum_{t=0}^{K-1}\left(Q(t) - \sum_{i=1}^N\sum_{j=1}^M p_{ij}^*(t)\right)^2\\
&+\sum_{i=1}^N\sum_{j=1}^M\sum_{t=0}^{K-1}c_{ij}\left(p_{ij}^*(t)-\overline{p_{ij}}(t)\right)^2,\label{opt_eq_2}
\end{flalign}
where $Q(t)$ represents the targeted aggregate power consumption level at time $t$. While the first term in \eqref{opt_eq_2} minimizes the deviation from $Q(t)$ at any time interval $t$, the second term minimizes the total social discomfort of the neighborhood. 

The coordination agent seeks a price vector minimizing its objective function. However, as the coordination agent has to act in accordance with some regulations, the electricity price cannot exceed a certain threshold. Similarly, the coordination agent cannot lower the price freely because of the cost of generating electricity. We assume that these requirements are satisfied as long as $\boldsymbol{\pi} \in \mathcal{S}$. The resulting bilevel optimization problem is:
\begin{align}
\min_{\pi \in \mathcal{S}}\quad &f(\boldsymbol{p}^*, \boldsymbol{\pi})\label{opt_eq_3}\\
\text{s.t.} \quad & \boldsymbol{p}_i^*:=\argmin_{\boldsymbol{p}_i\in \P_i}\ f_i(\boldsymbol{p}_i, \boldsymbol{\pi}), \quad \forall i.\label{opt_eq_5}
\end{align}

For a given price vector, the home level optimization problem is a convex quadratic program. Therefore, the bilevel problem given in Equations \eqref{opt_eq_3}-\eqref{opt_eq_5} can be reformulated as a single level problem by replacing the optimality constraint given in Equation \eqref{opt_eq_5} with the KKT conditions of each home:
\begin{align}
\min_{\substack{\boldsymbol{p}_1, \dots, \boldsymbol{p}_N,\\ \boldsymbol{\lambda}_1, \dots, \boldsymbol{\lambda}_N, \\\boldsymbol{\pi} \in \mathcal{S}}}\quad &f(\boldsymbol{p}, \boldsymbol{\pi})\label{opt_eq_6}\\[0.5em]
\text{s.t.} \qquad  &\nabla_{\boldsymbol{p}_i}f_i(\boldsymbol{p}_i, \bpi) + \bG_i'\boldsymbol{\lambda}_i = 0, \qquad \forall i,\label{opt_eq_7}\\
\qquad & \text{Diag}(\boldsymbol{\lambda}_i)\left(\bG_i\boldsymbol{p}_i-\bh_i\right) =0,\qquad \forall i,\label{opt_eq_8}\\
\qquad & \bG_i\boldsymbol{p}_i \leq h_i,\qquad\qquad\qquad\qquad\ \forall i,\label{opt_eq_9}\\ 
\qquad & \boldsymbol{\lambda}_i \geq 0,\qquad\qquad\qquad\qquad\qquad \forall i,\label{opt_eq_10}
\end{align}
where $\boldsymbol{\lambda}_i$ is a vector of the dual variables corresponding to the constraints in the set $\P_i$. Since the objective function of each home is quadratic, the stationarity constraints given in \eqref{opt_eq_7} are linear. However, the constraints in \eqref{opt_eq_8} enforcing the complementary slackness condition are non-convex. Thus, the optimization problem in Equations \eqref{opt_eq_6}-\eqref{opt_eq_10} is a non-convex {\em{Quadratically Constrained Quadratic Program}} (QCQP).  \cmmnt{which make the optimization problem given in Equations \eqref{opt_eq_6}-\eqref{opt_eq_10} hard to solve.}

\section{GRADIENT BASED OPTIMIZATION FRAMEWORK}

Although the commercial solvers can handle non-convex QCQP problems by transforming them into {\em{Mixed Integer Programs}} (MIP), solving MIP problems can be computationally intractable when the problem size grows. Therefore, we propose an efficient gradient-based optimization technique to solve the formulated problem. Since the coordination agent minimizes its objective function by choosing a price vector, the gradient that we need to compute is as follows:
\begin{flalign}
\nabla_{\bpi}f(\bp^*, \bpi)\cmmnt{\frac{df(\bp^*, \bpi)}{d\boldsymbol{\pi}}}=f_{\bpi}(\bp^*, \bpi) + f_{\bp}(\bp^*, \bpi)\nabla_{\bpi}\bp^*\cmmnt{'\frac{d\bp^*}{d\bpi}},
\end{flalign}
where $f_{\bpi}$ and $f_{\bp}$ represent the partial derivatives, and $\nabla_{\bpi}\bp^*$ is a Jacobian matrix showing how the neighborhood alters its optimal load schedule in response to any change in the price vector.\cmmnt{is a Jacobian matrix showing the change of the optimal response of the neighborhood to the price vector.} Since the objective function of the coordination agent does not involve the price vector, $f_{\bpi}$ is equal to $0$. Thus, the gradient calculation can be reduced to
\begin{flalign}
\label{eq_jacob}
\nabla_{\bpi}f(\bp^*, \bpi)\cmmnt{\frac{df(\bp^*, \bpi)}{d\boldsymbol{\pi}}}= f_{\bp}(\bp^*, \bpi)'\nabla_{\bpi}\bp^*.
\end{flalign}

While it is straightforward to calculate $f_{\bp}$, estimating $\nabla_{\bpi}\bp^*$ requires a bit more work. For a given $\bpi$, the optimization problem solved by home $i$, $f_i(\boldsymbol{p}_i, \boldsymbol{\pi})$, is a quadratic program with linear constraints. In \cite{amos2017optnet}, how the optimal solution responds to any change in the problem parameters (e.g., cost function, coefficients in constraint matrix) is shown by taking matrix differentials of the KKT conditions. Moreover, \cite{barratt2018differentiability} provides a detailed derivation of this approach by using the implicit function theorem. Finally, \cite{my_ref12} develops a gradient descent-based approach to optimize the leader's strategy in a Stackelberg game. Similar to \cite{my_ref12}, each home solves a convex optimization problem with a twice differentiable objective function in our case. Therefore, we describe how to obtain $\nabla_{\bpi}\bp^*$ by recalling the strategy illustrated in \cite{my_ref12}.

\cmmnt{\textbf{(We know that Slater' s condition is sufficient for strong duality in convex problems. The issue is that it is not necessary for strong duality to hold. So, we have to assume that Slater's condition holds to be able to use implicit function theorem as it is described. They published one paper in ICML 2017 another in AAAI 2022 without mentioning this important detail)}}Since the inequality constraints in Equation \eqref{opt_eq_9} are affine, the weak form of Slater's condition holds\cmmnt{see the last comment at: https://math.stackexchange.com/questions/3804265/weak-slaters-condition}. Thus, \cmmnt{If we assume that Slater's condition holds,} around $\bp_i^*$, the KKT conditions we have to maintain reduce to the stationarity and the complementary slackness conditions. Then, by taking the differentials of the remaining KKT conditions of home $i$, we obtain the following set of equations in matrix form:
\begin{equation}
\label{eq_KKT}
\begin{split}
 &\begin{bmatrix}
     \nabla^2_{\bp_i\bp_i}f_i & \bG_i'\\
     \text{Diag}(\boldsymbol{\lambda}_i^*)\bG_i & \text{Diag}(\bG_i\bp_i^* -\bh_i)\\ 
 \end{bmatrix}
 \begin{bmatrix}
     d\bp_i^*\\
     d\boldsymbol{\lambda_i}^* 
 \end{bmatrix} \\
  &=\begin{bmatrix}
  -\nabla^2_{\bpi\bp_i}f_i d\bpi - d\bG_i'\boldsymbol{\lambda}_i^*\\
  -\text{Diag}(\boldsymbol{\lambda}_i^*)(d\bG_i\bp_i^*-d\bh_i)
    \end{bmatrix}.
\end{split}
\end{equation}
If we set $d\bpi$ and $dG_i$ to the identity and the zero matrices of appropriate dimensions, respectively, and solve the system of equations above for $d\bp_i^*$, we obtain the Jacobian matrix  $\nabla_{\bpi}\bp_i^*$. However, as we want to obtain $\nabla_{\bpi}\bp^*$, we need to write the equations above for each home, and solve the following set of equations: 
\cmmnt{Solving the system of equations above for $d\bp_i^*$ after setting $d\bpi$ to the identity matrix of appropriate dimension, and $dG_i$ to zero in the right hand side of Equation \eqref{eq_KKT}, the Jacobian matrices showing the response of optimal primal and dual variables to any change in price vector can be constructed. } 
\begin{equation}
\label{eq_matrix}
 \begin{bmatrix}
     \boldsymbol{C} & \bG'\\
     \text{Diag}(\boldsymbol{\lambda}^*)\bG & \text{Diag}(\bG\bp -\bh)\\ 
 \end{bmatrix}
 \begin{bmatrix}
     d\bp^*\\
     d\boldsymbol{\lambda}^* 
 \end{bmatrix}
  = 
 \begin{bmatrix}
     -\nabla_{\bpi} \boldsymbol{v}\\
     0 
 \end{bmatrix},
\end{equation}
where $\boldsymbol{C} = \text{Diag}(\nabla^2_{\bp_1\bp_1}f_1,\ ...,\ \nabla^2_{\bp_N\bp_N}f_N)$, $\bG=\text{Diag}(\bG_1,\ ...,\ \bG_N)$, and $\boldsymbol{v} = \left[(\nabla_{\boldsymbol{p}_1}f_1)'\ ...\ (\nabla_{\boldsymbol{p}_N}f_N)'\right]'$. Unlike \cite{my_ref12}, since the optimal action of each follower does not depend on the other followers' actions in our problem, the large matrix on the left-hand side of Equation \eqref{eq_matrix} consists of diagonal blocks which do not interact with other followers' data. Hence, we can calculate $\nabla_{\bpi}\bp^*$ in a distributed way.

\begin{algorithm}[ht]
\caption{The framework of a distributed algorithm with a given maximum number of iterations $k^{\text{max}}$, and \cmmnt{threshold on gradient norm $\epsilon$.}threshold on objective improvement rate $\epsilon$.}
\label{alg:one}
\hspace*{\algorithmicindent}\textbf{Input:} Batch size $B$, learning rate $\alpha$, online optimization algorithm $\mathbb{A}$.
\begin{algorithmic}[1]
\State \textbf{Initialize}
\State Select a random price vector $\bpi_{1} \in \mathcal{S}$, and set $\bpi=\bpi_{1}$.
\State Set $z_0 = \infty$.
\For{\texttt{$k=1\ ...\ k^{\text{max}}$}}
        \State Set $\mathbb{B}=\emptyset$, and randomly select $B$ homes.  
        \State Store the indices of these homes in $\mathbb{B}$.
        \State Set $\boldsymbol{g}_k=0$.
        \ForEach{$i\ \in\ \mathbb{B}$}
            \State Solve Equation \eqref{opt_eq_home} to obtain $\bp_i^*$. 
            \State Solve Equation \eqref{eq_KKT} for $d\bp_i^*$ to obtain $\nabla_{\bpi}\bp_i^*$.
            \State $\boldsymbol{g}_k= \boldsymbol{g}_k + f_{\bp_i}(\bp^*, \bpi)' \nabla_{\bpi}\bp_i^*$.
        \EndFor
        \State Send $\boldsymbol{g}_k$, $\alpha$, and $\bpi$ to Algorithm $\mathbb{A}$.
        \State Receive $\bpi_{k+1}$ from Algorithm $\mathbb{A}$.
        \State $\bpi=$ Projection of $\bpi_{k+1}$ onto $\mathcal{S}$.
        \State $z_k = f(\bp^*,\bpi)$.
        \If {${|z_k - z_{k-1}|}/{z_{k-1}} \leq \epsilon$}
        \State \textbf{break}
        \EndIf
\EndFor
\State \Return $\bpi$.
\end{algorithmic}
\end{algorithm}
\cmmnt{which allows us to design a distributed optimization algorithm 
Hence, the large matrix on the left-hand side of Equation \eqref{eq_matrix} consists of diagonal blocks which do not interact with other homes' data, which allows us to design an optimization algorithm calculating $\frac{d\bp^*}{d\bpi}$ in a distributed way.}

In reality, we do not even need to construct the Jacobian matrix $\nabla_{\bpi}\bp^*$ since the gradient calculation in Equation \eqref{eq_jacob} is equivalent to:\cmmnt{we can rewrite the Equation \eqref{eq_jacob} in a slightly different way:}
\begin{flalign}
\nabla_{\bpi}f(\bp^*, \bpi)\cmmnt{\frac{df(\bp^*, \bpi)}{d\boldsymbol{\pi}}}&=
 \begin{bmatrix}
     f_{\bp_1}(\bp^*, \bpi)'\ ...\ f_{\bp_N}(\bp^*, \bpi)'
 \end{bmatrix}
  \begin{bmatrix}
  \nabla_{\bpi}\bp_1^* \\[-0.2em]
  \vdots \\[-0.2em]
  \nabla_{\bpi}\bp_N^*
  \end{bmatrix}\\[-0.6em]
 &=\sum_{i=1}^N f_{\bp_i}(\bp^*, \bpi)' \nabla_{\bpi}\bp_i^*,\label{eq_sgd}
\end{flalign}
where we can obtain $\nabla_{\bpi}\bp_i^*$ for each home $i$ by solving Equation \eqref{eq_KKT}  as it is already explained. The overall optimization problem we formulate is non-convex, and we know that the gradient based approaches can get stuck at local minima even if the step length is chosen appropriately. On the other hand, the Stochastic Gradient Descent (SGD)-based approaches may avoid these local minimums by utilizing the randomness during the global minimum search. Thus, they are extensively used to optimize non-convex loss functions, especially in studies involving artificial neural networks. Fortunately, we can employ SGD-based optimization techniques in our problem since it is possible to obtain unbiased estimators for $\nabla_{\bpi}f(\bp^*, \bpi)$\cmmnt{$\frac{df(\bp^*, \bpi)}{d\boldsymbol{\pi}}$}. We can make this observation quickly after manipulating Equation \eqref{eq_sgd} in the following way:
\begin{flalign}
\nabla_{\bpi}f(\bp^*, \bpi)\cmmnt{\frac{df(\bp^*, \bpi)}{d\boldsymbol{\pi}}}&= N\  \sum_{i=1}^N \frac{1}{N}f_{\bp_i}(\bp^*, \bpi)' \nabla_{\bpi}\bp_i^*\\
&=N\  \mathbb{E}_{i\sim z} \left[f_{\bp_i}(\bp^*, \bpi)' \nabla_{\bpi}\bp_i^*\right],
\end{flalign}
where $z$ denotes a Discrete Uniform$(1,N)$, and $\mathbb{E}$ is the expectation operator. 

The pseudocode showing the steps of the distributed optimization algorithm is given in Algorithm 1. Since each home independently performs the inner loop steps, we have a distributed optimization algorithm. Therefore, the designed algorithm solves the issues related to efficiency and data privacy.

\section{Experiments}

In this study, we propose not only a Stackelberg game that can keep the load consumption of the neighborhood around some target level without causing too much social discomfort but also a gradient based distributed optimization framework that can find effective solutions in a reasonable amount of time. In order to test these hypotheses, we conduct various experiments on a simulation environment.
First, the details of our experimental setup are explained in Section \ref{exp_sim}. Next, the optimization performance of the proposed optimization framework with different parameters is observed, and the optimization performance is compared with a commercial solver in Section \ref{exp_comp}. Finally, the capability of the proposed optimization problem to control the load consumption of the community is shown in Section \ref{exp_load}.

\subsection {Simulation Setup}
\label{exp_sim}

A simulation environment that can create a residential neighborhood, in which each home is equipped with the appliances mentioned in Section \ref{sec_app}, is designed in Python, and we utilize the Gurobi 9.5.1 package to solve both the centralized non-convex QCQP problem given in Equations \eqref{opt_eq_6}-\eqref{opt_eq_10} and the convex optimization problem given in\cmmnt{Equation} \eqref{opt_eq_home}. \cmmnt{we utilize the non-convex QCQP solver available in the Gurobi 9.5.1 package to solve the centralized optimization problem given in Equations \eqref{opt_eq_6}-\eqref{opt_eq_10}.}

\cmmnt{The properties of the appliances are identical across the neighborhood. However, both the preferences and the consumption habits of households can vary for each appliance, thus their desirable load schedules. Moreover, the penalty coefficient, $c_{ij}$, is randomly set for each home and appliance. Both the details of the appliance properties and the probability distributions that we use to sample the preferences of households are provided in \cite{this_paper_longer}. Additionally, the source code is \cmmnt{publicly} available\footnote{Code is publicly available at \url{https://github.com/erhancanozcan/energy-qcqp-grad}.}.}
The properties of the appliances are identical across the neighborhood. However, both the preferences and the consumption habits of households can vary for each appliance, thus their desirable load schedules. Moreover, the cost of deviation from the desirable load curve is different for each home and appliance. The probability distributions that we use to sample the preferences and penalty coefficients, $c_{ij}$s, of households as well as the details of the appliance properties are provided in \cite{this_paper_longer}. Additionally, the source code is publicly available\footnote{Code available at \url{https://github.com/erhancanozcan/energy-qcqp-grad}.}.

In our simulations, we assume that the length of the planning horizon is 24 hours, and each time interval lasts 15 minutes $(K=96)$. Since the desirable load consumption schedule of each home changes dynamically at each time interval, the optimization algorithm must provide a reasonably good price vector $\bpi$ in less than 15 minutes. We assume that the electricity price always has to be in a certain interval, and the price range we consider in our experiments is as follows:\[\mathcal{S}=\{\bpi\ |\ 0.1\leq \pi(t)\leq 1.0,\quad \forall t\}.\]\cmmnt{\[\Pi=\{\bpi\ |\ 0.1\leq \pi(t)\leq 1.0,\quad t=0,1,\cdots,K-1\}.\]}

Finally, $Q(t)$, the target load level the coordination agent wants to achieve at time interval $t$, is set as the time average of total desirable power consumption. In other words, $Q(t)$ takes the form \[Q(t)=\frac{\sum_{i=1}^N\sum_{j=1}^M\sum_{t=0}^{K-1}\overline{p_{ij}}(t)}{K},\ \forall t.\]

\subsection {Computational Performance of Distributed Optimization}
\label{exp_comp}

The proposed algorithm has three parameters, and the performance of the proposed algorithm changes with respect to these parameters. The first parameter is the batch size $B$, and we have two levels for this parameter in our experiments. It is either equal to $25$ or equal to the total number of homes in the community (e.g., full batch). Secondly, we consider the Adam optimizer and the scaled SGD (e.g., $\alpha\propto\frac{1}{\sqrt{k}}$, where $k$ is the iteration number) as two options for the online learning algorithm $\mathbb{A}$. Lastly, the possible $\alpha$ values we consider for each online learning algorithm vary. While two $\alpha$ values we consider for the scaled SGD are $1e^{-5}$ and $1e^{-6}$, the possible learning rate options for the Adam optimizer are $1e^{-1}$ and $1e^{0}$. In addition to these three parameters, we set $k^{\text{max}}$ to $50$, and $\epsilon$ to $1e^{-3}$ in our experiments. 

In order to observe how the performance of the proposed algorithm changes with respect to different parameter settings, we randomly initialize communities with varying number of homes (e.g., 50, 100, 250). To decrease the effect of random initialization, we use five different seeds in our experiments. Table 1 shows the number of times in which each algorithm yields the lowest objective value among all options across five trials.
\begin{table}[!h]
\centering
\caption{Lowest Objective Across Five Trials}
\label{table_one}
\renewcommand{\arraystretch}{.9}
\begin{tabular}{lccc}
\hline
\textbf{Algorithm 1 w/parameters}                & \textbf{50 homes} & \textbf{100 homes} & \textbf{250 homes}  \\ 
\hline
\textbf{Adam\_B25\_$\alpha$1en1}          & 2               & 2                & 3 \\

\textbf{Adam\_B25\_$\alpha$1en0}          & 0               & 0                & 0 \\

\textbf{Adam\_fullbatch\_$\alpha$1en1}          & 0               & 0                & 0 \\                          
\textbf{Adam\_fullbatch\_$\alpha$1en0}          & 0               & 0                & 0 \\ 

\textbf{scaledSGD\_B25\_$\alpha$1en5}          & 1               & 1                & 1 \\

\textbf{scaledSGD\_B25\_$\alpha$1en6}          & 0               & 1                & 1 \\

\textbf{scaledSGD\_fullbatch\_$\alpha$1en5}          & 0               & 0                & 0 \\                          
\textbf{scaledSGD\_fullbatch\_$\alpha$1en6}          & 2               & 1                & 0 \\ 

\hline
\end{tabular}
\end{table}

\begin{figure}[H]
\centering
\includegraphics[width=3.5in]{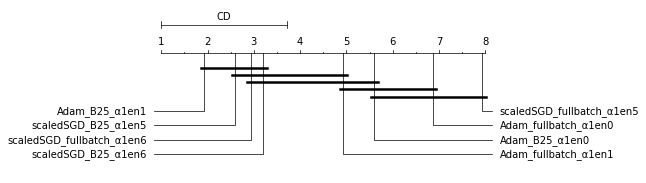}
\caption{The average ranks over attained objective value for all competing methods on 15 problems. Lower rank corresponds to a better performance. Bold horizontal lines shows the methods that are not different in a 5\% significance level according to the Nemenyi test.}
\label{fig_rank}
\end{figure}

According to Table \ref {table_one}, a mini batch optimization technique helps our algorithm to attain better objective values. On the other hand, both the Adam optimizer and the scaled SGD may yield the lowest objective value depending on the problem at hand. In order to decide which setup yields the best performance in Algorithm \ref{alg:one}, we compare the competing setups based on their average ranks over attained objective value. Figure \ref{fig_rank} provides the average ranks of each competing setting with critical differences where lower rank corresponds to a better performance.

According to Figure \ref{fig_rank}, the Adam optimizer with mini batch seems as a good alternative to attain lower objective values across a range of problems although there is no statistically significant difference with the setups using the scaled SGD as the optimization algorithm. We pick the top three options in the average rank plot as our candidate setups, and continue our analysis with those. 

It is possible to solve the centralized version of the formulated problem by using the non-convex QCQP solver available in Gurobi, and this strategy is a natural competitor of our algorithm. However, as the centralized version is an MIQP problem, it is hard to solve this problem optimally. Unfortunately, Gurobi fails to solve this problem, and it cannot even find a feasible point within the desired time limit in many trials. The optimization iterations of our algorithm start using a random price vector, and a corresponding optimal load schedule plan of the community. By using this information, an initial feasible solution can be easily found for the MIQP problem, and we provide this information to Gurobi in order to have a fair comparison between our algorithm and Gurobi. Even in this case, we observe either a slight improvement in the objective or no improvement at all when Gurobi is employed. Suppose $z^\text{G}$ denotes the objective value calculated by Gurobi at the end of 900 seconds for a given problem, and the best objective value that the k-th method can find within this time period is denoted by $z^{\text{method}_k}$. Then, we define the improvement ratio of the k-th method as follows: \[\text{Ratio (}\text{method}_k) = \frac{z^{\text{G}}-z^{\text{method}_k}}{z^{\text{method}_k}},\ \forall\ k.\]
Table \ref{table_three} provides the mean improvement ratio for each algorithm.

\begin{table}[!h]
\centering
\caption{Average Improvement Ratio Compared to Gurobi}
\label{table_three}
\renewcommand{\arraystretch}{.9}
\begin{tabular}{lccc}
\hline
\textbf{Algorithm 1 w/parameters}                & \textbf{50 homes} & \textbf{100 homes} & \textbf{250 homes}  \\ 
\hline
\textbf{Adam\_B25\_$\alpha$1en1}          & 35.2               & 35.9                & 34.6 \\



\textbf{scaledSGD\_B25\_$\alpha$1en5}          & 32.7               & 33.1                & 34.4 \\


\textbf{scaledSGD\_fullbatch\_$\alpha$1en6}          & 34.9               & 35.6                & 30.4 \\ 

\hline
\end{tabular}
\end{table}

According to Table \ref{table_three}, on average, the objective values provided by our algorithm are at least 30 times better than the objective values provided by Gurobi. Moreover, the run time of our algorithm is significantly shorter than Gurobi. Table \ref{table_two} summarizes the average run time of our algorithm under different parameter settings. In all cases, our algorithm can provide the solution within the time limit.
\cmmnt{
Therefore, we compare the quality of the solutions provided by our algorithm and Gurobi based on a pseudo optimality gap. At each iteration, Gurobi provides a lower bound for the problem by using branch and bound technique, and suppose that $z^\text{LB}$ denotes the lower bound calculated by Gurobi at the end of 900 seconds for a given problem. On the other hand, the best objective value that the k-th method can find within this time period is denoted by $z^{\text{method}_k}$. Then, we define the pseudogap of the k-th method as follows: \[\text{pseudogap (}\text{method}_k) = \frac{z^{\text{method}_k}-z^\text{LB}}{z^\text{LB}},\ \forall\ k.\]
Figure ....... shows the mean pseudogap (along with a 95\% confidence interval) for each algorithm.

Based on Figure ..., our algorithm can provide solutions with significantly small pseudogaps compared to Gurobi. Hence, our approach provides more effective solutions than Gurobi. Finally, Table \ref{table_two} summarizes the average run time of our algorithm under different parameter settings. In all cases, our algorithm can provide the solution within the time limit.} 
\begin{table}[!h]
\centering
\caption{Average Run Time in Seconds}
\label{table_two}
\renewcommand{\arraystretch}{.9}
\begin{tabular}{lccc}
\hline
\textbf{Algorithm 1 w/parameters}                & \textbf{50 homes} & \textbf{100 homes} & \textbf{250 homes}  \\ 
\hline
\textbf{Adam\_B25\_$\alpha$1en1}          & 130.9               & 165.3                & 225.2 \\



\textbf{scaledSGD\_B25\_$\alpha$1en5}          & 95.8               & 116.9                & 222.4 \\


\textbf{scaledSGD\_fullbatch\_$\alpha$1en6}          & 65.7               & 41.1                & 73.3 \\ 

\hline
\end{tabular}
\end{table}

\subsection {The Load Shaping Capacity}
\label{exp_load}

\begin{figure}[H]
\centering
\includegraphics[width=3.25in]{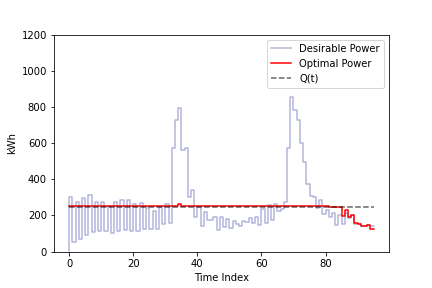}
\caption{The desirable and the optimal power consumption profiles of the community consisting of 100 homes.}
\label{fig_pp}
\end{figure}

In this section, we show the capability of the proposed formulation to shape the overall load consumption of the community after optimizing the price vector via the proposed algorithm. Initially, the desirable load consumption of home $i$ at time $t$ is equal to $\sum_{j=1}^M\overline{p_{ij}}(t)$, and the desirable load consumption of the whole community is equal to $\sum_{i=1}^N\sum_{j=1}^M\overline{p_{ij}}(t).$ On the other hand, $\sum_{i=1}^N\sum_{j=1}^M p_{ij}^*(t)$ denotes the optimal load consumption of the community at time $t$ with respect to the optimized price vector $\bpi^*$. Figure \ref{fig_pp} shows how the overall desirable load consumption of the community changes when the coordination agent announces an optimized price vector.
According to Figure \ref{fig_pp}, the observed peaks in the desirable load consumption profile of community are prevented by setting the electricity price vector appropriately. In Figure \ref{fig_pp}, the horizontal dashed line denotes the targeted aggregate power consumption level set by the coordination agent, and we observe that the optimal load consumption curve of the community mostly matches this target. Therefore, our approach can be helpful to maintain the balance between electricity demand and supply as well as mitigating the peak load problem.

\section{CONCLUSIONS}

In this paper, we present a Stackelberg game to control the overall load consumption of a residential neighborhood. As the control decisions are taken by considering the individual preferences of the households, we ensure that the user comfort is maintained at any time.\cmmnt{As the control decisions are taken by individual homes inline with their preferences, we ensure that the user comfort is maintained at any time.} However, refining the appliance models to reflect the user comfort better is one of our future goals. Although it can be challenging to solve the formulated problem, the proposed distributed optimization framework can obtain effective solutions in a significantly shorter time period. As future work, we want to investigate how we can use our algorithm to help Gurobi to start its optimization from a better branch, as part of its branch and bound search.

\bibliographystyle{IEEEtran} 
\bibliography{IEEEabrv,IEEEexample}

\appendix

The set $\P_i$ consists of the comfort related constraints of home $i$, and the constraints for each appliance are introduced in this section. The appliance models we consider in this study are very similar to the appliances in \cite{prev_paper}. While the properties of the appliances and the probability distributions that we use to sample the desirable load schedules for homes are the same, there are minor differences in the constraints. Thus, we use a similar notation to discuss those appliances whenever it is possible. 

The first modification that we have to apply is to simplify the work mechanism of ON/OFF type appliances by removing binary variables. To be able to use the proposed gradient optimization scheme, we have to transform the bilevel optimization problem into single level problem by using KKT conditions. Therefore, the home level optimization problem has to be a convex problem, which is not the case when binary variables are utilized. Secondly, our approach assumes that the comfort related constraints are defined by using linear inequalities. Therefore, we modify the constraints given in \cite{prev_paper} by using traditional linear optimization techniques. 

In our study, the planning horizon is divided into $K$ equal time intervals, and the whole planning horizon is spanned by the set $\T$, where \[\T:=\{0,1,\dots,K-1\}.\]

\subsection{Heating, Ventilation, and Air Conditioning}

The comfortable thermal range for the consumer $i$ is between $T_i^{low}$ and  $T_i^{upper}$. Therefore, the HVAC system has to ensure that the room temperature remains in this interval at any time. The HVAC system has heating and cooling modes, and the HVAC activity depends on this mode as well as the room temperature. For example, the HVAC gets activated when it is in its heating mode and the room temperature is less than $T_i^{low}$. On the contrary, the HVAC starts working if the room temperature is greater than $T_i^{upper}$ and it is in its cooling mode. 

Suppose that $T_{i,in}(t)$ denotes the room temperature in home $i$ at time $t$, and $T_{out}(t)$ is the outside temperature at time $t$. When the HVAC is in its heating mode, the temperature update rule used in \cite{prev_paper} is as follows:
\begin{equation}
\begin{split}
\label{hvac_recursive}
T_{i,in}(t+1)=&T_{i,in}(t) +\gamma_1(T_{out}(t)-T_{i,in} (t))\\ 
&+\gamma_2p_{i,\text{HVAC}}(t),\quad\quad\qquad t\in \T,
\end{split}
\end{equation}
where $\gamma_1,\gamma_2$ are the home specific coefficients related to building isolation and heating or cooling gain, respectively. In this representation, the room temperature is calculated recursively. However, given the initial room temperature, $T_{i,in} (0)$, the room temperature at time $t$ can also be expressed in closed form by expanding the recursive equation:
\begin{equation}
\begin{split}
\label{hvac_closed}
T_{i,in}(t)&=(1-\gamma_1)^tT_{i,in}(0)\\
&+\sum_{a=0}^{t-1}(1-\gamma_1)^a\gamma_1T_{out}(t-1-a)\\
&+\sum_{a=0}^{t-1}(1-\gamma_1)^a\gamma_2 p_{i,HVAC}(t-1-a),\ \forall t\in\T.
\end{split}
\end{equation}

Suppose that the nominal power consumption of the HVAC is equal to $N_{\text{HVAC}}$ whenever it gets activated. However, modeling this behaviour requires using binary variables. Instead, we relax this assumption by forcing the load consumption of the HVAC to be less than $N_{\text{HVAC}}$. We believe this relaxation is reasonable since $\frac{p_{i,\text{HVAC}}(t)}{N_{\text{HVAC}}}$ can be interpreted as the percentage of time on which the HVAC remains in the ON status at time interval $t$. Then, HEMS can ensure the thermal comfort of the user $i$ by considering the following constraints:
\begin{align}
&T_i^{low}\leq T_{i,in}(t) \leq T_i^{upper},\quad t\in \T,\label{hvac_a} \\
\quad& p_{i,\text{HVAC}}(t) \leq N_{\text{HVAC}},\quad t\in \T\label{hvac_d},\\
\quad& p_{i,\text{HVAC}}(t) \geq 0,\quad t\in \T\label{hvac_e},
\end{align}
\cmmnt{
Suppose that the nominal power consumption of the HVAC is equal to $N_{\text{HVAC}}$ whenever it gets activated. Then, $\frac{p_{i,\text{HVAC}}(t)}{N_{\text{HVAC}}}$ denotes the ON/OFF status of the HVAC at time t. Finally, suppose that $T_i^{set}(t)$ denotes the set temperature that the user has select to operate the HVAC. Then, HEMS can ensure the thermal comfort of the user $i$ by considering the following constraints:
\begin{align}
&T_i^{low}\leq T_{i,in}(t) \leq T_i^{upper},\quad t\in \T,\label{hvac_a} \\
& T_i^{set}(t) \geq T_{i,in}(t) - M^{H}\left(1-\frac{p_{i,\text{HVAC}}(t)}{N_{\text{HVAC}}}\right), \quad t\in \T,\label{hvac_b}\\
& T_{i,in}(t) \geq T_i^{set}(t)  - M^{H} \frac{p_{i,\text{HVAC}}(t)}{N_{\text{HVAC}}},\quad t\in \T,\label{hvac_c}\\
\quad& p_{i,\text{HVAC}}(t) \leq N_{\text{HVAC}},\quad t\in \T\label{hvac_d},\\
\quad& p_{i,\text{HVAC}}(t) \geq 0,\quad t\in \T\label{hvac_e},
\end{align}}
where Equation \eqref{hvac_a} ensures that the room temperature remains in the comfortable thermal range at time $t$, and Equations \eqref{hvac_d} specifies the maximum load consumption of the HVAC. \cmmnt{Equations \eqref{hvac_b},  \eqref{hvac_c}, and \eqref{hvac_d} relate the set temperature and ON/OFF status of the HVAC with the help of some large number $M^H$.}

While using the HVAC in its cooling mode, we need to flip the sign of the last term in the right-hand side of Equation \eqref{hvac_closed} since the HVAC is responsible from cooling the room. Hence, in the cooling mode, the room temperature is calculated by using the following closed form expression: 
\begin{equation}
\begin{split}
T_{i,in}(t)&=(1-\gamma_1)^tT_{i,in}(0)\\
&+\sum_{a=0}^{t-1}(1-\gamma_1)^a\gamma_1T_{out}(t-1-a)\\
&-\sum_{a=0}^{t-1}(1-\gamma_1)^a\gamma_2 p_{i,HVAC}(t-1-a),\ \forall t\in\T.
\end{split}
\end{equation}
\cmmnt{Finally, the constraints in Equations \eqref{hvac_b} and \eqref{hvac_c} must be replaced by the following two equations:
\begin{align}
&T_{i,in}(t) \geq T_i^{set}(t) -M^{H}\left(1-\frac{p_{i,\text{HVAC}}(t)}{N_{\text{HVAC}}}\right),\\
&T_{i,in}(t) \leq T_i^{set}(t)  +M^{H}\frac{p_{i,\text{HVAC}}(t)}{N_{\text{HVAC}}}.
\end{align}
}
\subsection{Electric Water Heater (EWH)}

The role of HEMS is to ensure that there exists enough hot water in the water tank to satisfy the hot water demand of the user $i$ at any time point $t$, which is denoted by $y_{i,EWH}(t)$. Suppose that $C_{i,\text{EWH}}$, and $p_{\text{EHW}}^{max_i}$ represent the capacity of the tank and the maximum power the EWH can consume, respectively. Finally, we define the following decision variables:
\begin{itemize}
    \item $x_{i,\text{EWH}}(t)$ denotes the amount of hot water available in tank at time $t$.
    \item $z_{i,\text{EWH}}(t)$ denotes the amount of water that the EWH will heat from time $t$ to $(t+1)$
\end{itemize}
We use the relation between between heat transfer and temperature change to calculate the amount of water that the EWH can heat from time $t$ to $(t+1)$. Then, $z_{i,\text{EWH}}(t)$ can be expressed as follows:
\begin{equation}
\label{ewh_heating}
z_{i,\text{EWH}}(t) = \frac{p_{i,\text{EWH}}(t)\; \eta_{\text{EWH}}}{\rho(T_{d} - T_{t})},
\end{equation}
where $\eta_{\text{EWH}}$ is the water heating efficiency of EWH, $\rho$ is the specific heat of water\cmmnt{ in $\left(J / kg^\circ C\right)$}, and $T_{d}$ is the desired water temperature\cmmnt{ in $(^\circ C)$}, and $T_{t}$ is the tap water temperature\cmmnt{ in $(^\circ C)$}.

According to these definitions, we can write the following recursive expression to express the available hot water at time point $t$: 
\begin{equation}
\begin{split}
&x_{i,\text{EWH}}(t) = x_{i,\text{EWH}}(t-1) + z_{i,\text{EWH}}(t-1)\\ 
&\ \quad\qquad\qquad\qquad\qquad - y_{i,EWH}(t-1),  t \in \T.\label{ewh_recursive}
\end{split}
\end{equation}
As in the HVAC, we can expand the recursive function given in Equation \eqref{ewh_recursive}, and express the available hot water at time point $t$ as follows:
\begin{equation}
\begin{split}
x_{i,\text{EWH}}(t) =& x_{i,\text{EWH}}(0)\\ 
&+ \sum_{a=0}^{t-1} \left(z_{i,\text{EWH}}(a) - y_{i,EWH}(a)\right),  t \in \T,
\end{split}
\end{equation}
where $z_{i,\text{EWH}}(a)$ is calculated as in Equation \eqref{ewh_heating}.

Based on these definitions, we write the following set of constraints to maintain the user comfort:
\begin{align}
&x_{i,\text{EWH}}(t) \leq C_{i,\text{EWH}}, \quad t \in \T, \label{ewh_a}\\
&x_{i,\text{EWH}}(t) \geq y_{i,\text{EWH}}(t), \quad t \in \T,\label{ewh_b}\\
&p_{i,\text{EWH}}(t) \leq p_{\text{EHW}}^{max_i},\quad t \in \T,\label{ewh_c}\\
&p_{i,\text{EWH}}(t) \geq 0,\quad t \in \T,\label{ewh_d}\\
&x_{i,\text{EWH}}(t) \geq 0,\quad t \in \T,\label{ewh_f}
\end{align}
where Equation \eqref{ewh_a} ensures that the available hot water is less than the capacity of the tank, Equation \eqref{ewh_b} ensures that the hot water demand of the user is satisfied, Equations \eqref{ewh_c} and \eqref{ewh_d} set the lower and upper bounds for the power consumption of the EWH.

\subsection{Electric Vehicle (EV)}

Similar to EWH, the role of HEMS is to ensure that there exists enough charge in the battery of the EV to satisfy the demand of the user $i$ at any time point $t$, which is denoted by $y_{i,EV}(t)$. While the power capacity of the EV battery is denoted by $C_{i,\text{EV}}$, the maximum power that can be delivered to the battery at time $t$ is represented by $p_{\text{EV}}^{max_i}$. Finally, we define the following decision variables for the EV:
\begin{itemize}
    \item $x_{i,\text{EV}}(t)$ denotes the amount of available power in the battery at time $t$.
    \item $z_{i,\text{EV}}(t)$ denotes the amount of current delivered to the battery of the EV from time $t$ to $(t+1)$.
\end{itemize}

Unlike EWH, the battery of the car cannot be charged when it is in use. Hence, $z_{i,\text{EV}}(t)$ is defined as follows:
\[ z_{i,\text{EV}}(t)=\begin{cases} 
      p_{i,\text{EV}}(t), & \text{if } t\notin \T^{\text{EV}_i} \\
      0, & \text{otherwise}
   \end{cases},
\]
where $\T^{\text{EV}_i}:=\{t\ |\ y_{i,\text{EV}}(t)>0,\ t\in \T\}$. Based on these definitions, the amount of available power in the battery can be expressed as follows:
\begin{equation}
\begin{split}
x_{i,\text{EV}}(t) =& x_{i,\text{EV}}(0)\\ 
&+ \sum_{a=0}^{t-1} \left(z_{i,\text{EV}}(a) - y_{i,EV}(a)\right),  t \in \T.
\end{split}
\end{equation}
Then, the user comfort can be maintained by writing the the following set of constraints:
\begin{align}
& x_{i,\text{EV}}(t) \leq C_{i,\text{EV}}, \quad t \in \T, \label{ev_a}\\
& x_{i,\text{EV}}(t) \geq y_{i,\text{EV}}(t), \quad t \in \T, \label{ev_b}\\
& p_{i,\text{EV}}(t) \leq p_{\text{EV}}^{max_i}, \quad t \in \T, \label{ev_c}\\
& p_{i,\text{EV}}(t) \geq 0, \quad t \in \T, \label{ev_d}
\end{align}
where Equation \eqref{ev_a} ensures that the amount of power stored in the battery is less than the capacity, Equation \eqref{ev_b} forces that the EV battery has enough charge to satisfy the customer needs, Equation \eqref{ev_c} and \eqref{ev_d} set the limits for the power consumption of the EV.

\subsection{Basic Appliances}

We consider basic appliances such as washing machine (WM), oven, and dryer in this study. In this section, we only introduce the constraint set of the WM, but these constraints can be easily generalized for other appliances. The study in \cite{prev_paper} uses binary variables to model the ON/OFF status of the corresponding appliance. Thus, the overall set of constraints ensuring the user comfort of the homes is non-convex. In this paper, the home level optimization problem is required to be a convex problem as we use KKT conditions to reformulate the bilevel problem. Therefore, we use simpler constraints to maintain the consumer comfort for the basic appliances.
Suppose that home $i$ wants WM operation to be completed between $\alpha^{start}_{i,\text{WM}}$ and $\alpha^{end}_{i,\text{WM}}$, and $C_{i,\text{WM}}$ denotes the total amount of power the WM consumes during the operation. Finally, $N_{i,\text{WM}}$ is the maximum amount of power the WM can consume at time interval $t$. Then, the user comfort can be maintained by writing the the following set of constraints:

\begin{align}
& \sum_{t \in \T^{\text{WM}_i}} p_{i,\text{WM}}(t) \leq C_{i,\text{WM}}, \label{wm_a}\\
& \sum_{t \in \T^{\text{WM}_i}} p_{i,\text{WM}}(t) \geq C_{i,\text{WM}}, \label{wm_b}\\
& p_{i,\text{WM}}(t) \leq N_{i,\text{WM}},\quad t\in \T,\label{wm_c}\\
& p_{i,\text{WM}}(t) \geq 0,\quad\quad\quad t\in \T,\label{wm_d}
\end{align}
where $\T^{\text{WM}_i}:=\{t\ |\ \alpha^{start}_{i,\text{WM}}\leq t \leq \alpha^{end}_{i,\text{WM}},\ t\in \T\}$. Equations \eqref{wm_a} and \eqref{wm_b} together imply that the amount of power consumed by the WM is equal to $C_{i,\text{WM}}$, and Equation \eqref{wm_c} limits the power consumption of the WM at time interval $t$.

\end{document}